\documentclass[12pt]{article}
\usepackage{amssymb}
\usepackage{amscd}
\usepackage{amsmath}
\usepackage{amsthm}
\usepackage{verbatim}

\let\geq\geqslant
\let\leq\leqslant

\newcommand{\g}{{\gamma}}
\renewcommand{\d}{{\delta}}
\newcommand{\e}{{\varepsilon}}

\newcommand{\f}{{\varphi}}

\newcommand{\s}{{\sigma}}

\newcommand{\G}{{\Gamma}}

\DeclareMathOperator{\sign}{sign}

\DeclareMathOperator{\Ran}{Ran}

\DeclareMathOperator{\Ker}{Ker}
\DeclareMathOperator{\Tr}{Tr}

\renewcommand\Im{\hbox{{\rm Im}}\,}
\renewcommand\Re{\hbox{{\rm Re}}\,}
\newcommand{\abs}[1]{\lvert#1\rvert}

\newcommand{\norm}[1]{\lVert#1\rVert}
\newcommand{\nnorm}[1]{\left\lVert#1\right\rVert}

\newcommand{\R}{{\mathbb R}}

\renewcommand{\H}{{\mathcal H}}
\newcommand{\K}{{\mathcal K}}

\numberwithin{equation}{section}

\theoremstyle{plain}
\newtheorem{theorem}{\bf Theorem}
\newtheorem{lemma}[theorem]{\bf Lemma}
\newtheorem{proposition}[theorem]{\bf Proposition}

\newtheorem{hypothesis}{\bf Hypothesis}

\theoremstyle{definition}

\theoremstyle{remark}
\newtheorem*{remark*}{\bf Remark}

\newcommand{\h}{\mathfrak{h}}
\newcommand{\wt}{\widetilde}
\newcommand{\F}{\mathcal{F}}

\begin{document}

\title{Differences of spectral projections and scattering matrix}
\author{A.~Pushnitski
\thanks{Department of Mathematics,
King's College London, 
Strand, London, WC2R  2LS, U.K.
email: \mbox{alexander.pushnitski@kcl.ac.uk}}
}
\maketitle


\begin{abstract}
In the scattering theory framework, we point out a connection between the 
spectrum of the scattering matrix of two operators and the spectrum of the
difference of spectral projections of these operators. 
\end{abstract}

\section{Introduction and results}

\textbf{1. Motivation and an informal description of results. }
Let $H_0$ and $H$ be self-adjoint operators in a Hilbert space $\H$ and suppose that the difference 
$V=H-H_0$ is a compact operator.  
For $\lambda\in\R$, we denote by $E_0(\lambda)$ and $E(\lambda)$ the spectral projections
of $H_0$ and $H$, corresponding to the interval  $(-\infty,\lambda)$. 
Our aim is to discuss the spectral properties 
of the operators 
\begin{equation}
D(\lambda)=E(\lambda)-E_0(\lambda),
\quad \lambda\in\R
\label{a1}
\end{equation}
and to point out the connection between these properties and the scattering matrix $S(\lambda)$
for the pair of operators $H_0$, $H$. 

It is well known that due to the compactness of $V$, for any continuous function
$\f$ which tends to zero at infinity, the difference 
\begin{equation}
\f(H)-\f(H_0)
\label{a2}
\end{equation}
is compact. However, the difference \eqref{a2} in general fails to be compact if 
$\f$ has discontinuities on the essential spectrum of $H_0$ and $H$.  
This observation goes back to M.~G.~Krein \cite{Krein} and was recently revisited in \cite{KM}; 
we will say more on this in section 1.3. An attempt to understand Krein's example was part
of the motivation for this paper. 

The first question we address is the nature of the essential spectrum of the operators
$D(\lambda)$, as these are the simplest operators of the type \eqref{a2} when $\f$ has a 
discontinuity. We consider this problem in the scattering theory framework, i.e. we make
certain typical for the scattering theory assumptions of the Kato smoothness type. 
These assumptions, in particular,
ensure that the scattering matrix $S(\lambda)$ for the pair $H_0$, $H$ is well defined. 
Under these assumptions, we prove (see Theorem~\ref{th1}) that 
\begin{equation}
\s_{ess}(D(\lambda))=[-a,a], 
\quad
a=\frac12\norm{S(\lambda)-I_\lambda}.
\label{a3}
\end{equation}
Here the scattering matrix $S(\lambda)$ acts in the fiber Hilbert space $\h(\lambda)$, which appears
in the diagonalisation of the absolutely continuous part of $H_0$ (see \eqref{a12} below) and
$I_\lambda$ is the identity operator in $\h(\lambda)$.
In particular, \eqref{a3} says that $D(\lambda)$ is compact if and only if $S(\lambda)=I_\lambda$.

Next, we consider the difference $D(\lambda)$ in the framework of the trace class scattering 
theory. Assuming that a certain  trace class condition on $V$ is fulfilled,  
we describe the a.c. spectrum of the operator $D(\lambda)$
in terms of the spectrum of the scattering matrix. 
See Theorem~\ref{th2} for the precise statement. 

Note that the question of the spectral analysis of the difference $D(\lambda)$ is well posed
regardless of any scattering theory type assumptions on the pair of operators $H_0$, $H$. 
Thus, the observations presented here might offer an insight into possible extensions of 
some elements of the scattering theory framework to wider classes of pairs of operators.

In this paper, we do not aim to prove our results under the optimal assumptions on $H_0$ and $H$.
Our aim is rather to point out the connection between 
the spectral properties of $D(\lambda)$ and $S(\lambda)$ while keeping the technical 
details simple.

Our construction borrows several ideas from the spectral theory of Hankel 
operators; see \cite{Power,Howland1,Howland2,Howland3}.

We denote by ${\mathfrak S}_\infty$ the class of all compact operators and by 
${\mathfrak S}_1$ and ${\mathfrak S}_2$ the trace class and the Hilbert-Schmidt class respectively. 
Along with the notation $E_0(\lambda)$, $E(\lambda)$ for $\lambda\in \R$, 
we also use the notation $E_0(\delta)$, $E(\delta)$ for the spectral projections 
of $H_0$ and $H$ associated with a Borel set $\delta\subset\R$.

\textbf{2. Statement of Results. }
Let $H-H_0=V=G^*V_0 G$, where $G$ is a bounded operator from $\H$ to an auxiliary  Hilbert space $\K$,
and $V_0$ is a bounded self-adjoint operator in $\K$.
The simplest case of such a  factorisation is when $\K=\H$, 
$G=\abs{V}^{1/2}$ and $V_0=\sign(V)$. 
Let us define 
\begin{equation}
F_0(\lambda)=GE_0(\lambda)G^*,
\quad
F(\lambda)=GE(\lambda)G^*,
\quad \lambda\in\R.
\label{a4}
\end{equation}
Next, let $\delta\subset\sigma_{ac}(H_0)$ be an open interval. 
\begin{hypothesis}
\label{hyp1}
The operator $G$ is compact. 
For all $\lambda\in\delta$, the derivatives 
$F'_0(\lambda)=\frac{d}{d\lambda}F_0(\lambda)$ and $F'(\lambda)=\frac{d}{d\lambda}F(\lambda)$
exist in operator norm. The maps $\delta\ni \lambda \mapsto F_0'(\lambda)$ 
and $\delta\ni \lambda \mapsto F'(\lambda)$ are H\"older continuous (with some positive exponent)
in the operator norm. 
\end{hypothesis}
Hypothesis~\ref{hyp1} is close to (but stronger than) the local Kato smoothness assumption in scattering theory
(see \cite{RS3} or \cite{Yafaev}).
In fact, one can make the required assumption concerning $F'_0(\lambda)$ and in addition
assume that 
\begin{equation}
\lim_{\epsilon\to+0}(I+V_0 G(H_0-\lambda-i\epsilon)^{-1}G^*)
\text{ is invertible for all $\lambda\in \delta$.}
\label{aa}
\end{equation}
This will ensure that the required assumption holds true also for $F'(\lambda)$.

Next, we recall the definition of the scattering matrix. 
Let $\H^{(ac)}_0(\delta)\subset \Ran E_0(\delta)$ 
be the absolutely continuous subspace of the operator $H_0\mid E_0(\delta)$ and
$\H^{(ac)}(\delta)$ be the absolutely continuous subspace of  $H\mid E(\delta)$; let 
$P_0^{(ac)}$ be the orthogonal projection onto $\H^{(ac)}_0(\delta)$ in $\H$. 
Hypothesis~\ref{hyp1} ensures that the local wave operators 
$$
W_\pm :=s-\lim_{t\to \pm\infty} e^{itH}e^{-itH_0}P_0^{(ac)}
$$
exist and are complete: $\Ran W_\pm =\H^{(ac)}(\delta)$. The local scattering operator 
$S=W_+^*W_-$ is unitary in $\H^{(ac)}_0(\delta)$ and commutes with $H_0\mid\H^{(ac)}_0(\delta)$. 
Consider the direct integral decomposition  
\begin{equation}
\H^{(ac)}_0(\d)=\int_\delta^\oplus \h(\lambda)d\lambda
\label{a12}
\end{equation}
which diagonalises $H_0\mid\H^{(ac)}_0(\delta)$.
Then 
$$
S=\int_\delta^\oplus  S(\lambda)d\lambda, 
\quad 
S(\lambda): \h(\lambda)\to \h(\lambda). 
$$
The scattering matrix $S(\lambda)$ is unitary in $\h(\lambda)$. 
The compactness of $G$ ensures  
that $S(\lambda)-I_\lambda$ is compact for all $\lambda\in\delta$.

\begin{theorem}\label{th1}
Suppose that for some open interval $\delta\subset\R$, Hypothesis~\ref{hyp1}
holds true.  
Then for all $\lambda\in\delta$ formula \eqref{a3} holds true. 
\end{theorem}

Next, we describe the trace class result. 
Instead of Hypothesis~\ref{hyp1}, we need the following stronger hypothesis: 
\begin{hypothesis}
\label{hyp2}
The operator $G$ is Hilbert-Schmidt. 
For all $\lambda\in\delta$, the derivatives 
$F'_0(\lambda)$ and $F'(\lambda)$
exist in the trace norm. The maps $\delta\ni \lambda \mapsto F_0'(\lambda)$ 
and $\delta\ni \lambda \mapsto F'(\lambda)$ are H\"older continuous (with some positive exponent)
in the trace norm. 
\end{hypothesis}
Again, it suffices to assume the existence and  H\"older continuity of $F'_0$ and \eqref{aa}; then
$F'$ also exists and is H\"older continuous. 

Under Hypothesis~\ref{hyp2}, the operator $S(\lambda)-I_\lambda$ is compact for all $\lambda\in \delta$. 
Thus, the spectrum of $S(\lambda)$ consists of eigenvalues on the unit circle which can only accumulate
to $1$. 
For $\lambda\in\delta$, let $e^{i\theta_n(\lambda)}$, $\theta_n(\lambda)\in(0,2\pi)$, 
 be the eigenvalues of $S(\lambda)$ distinct from $1$. There may be finitely or infinitely many 
 of these eigenvalues.

\begin{theorem}\label{th2}
Suppose that for an open interval $\delta\subset\R$, 
Hypothesis~\ref{hyp2} holds true. 
Then for all $\lambda\in\delta$ the a.c. part of the operator $D(\lambda)$
is unitarily equivalent to a direct sum of operators of multiplication by $x$
in $L^2([-a_n,a_n],dx)$, $a_n=\frac12\abs{e^{i\theta_n(\lambda)}-1}=\sin(\theta_n(\lambda)/2)$. 
\end{theorem}

Using Theorems~\ref{th1} and \ref{th2}, one can also analyse the spectra of the operators 
\eqref{a2} for certain classes of piecewise continuous functions $\f$.

\textbf{3. Krein's Example.} 
In \cite{Krein}, M.~G.~Krein considers an example of the operator $H_0$ 
in $L^2(0,\infty)$ with the integral kernel $H_0(x,y)$ given by 
$$
H_0(x,y)=
\begin{cases}
\sinh(x) e^{-y}, \quad x\leq y,
\\
\sinh(y) e^{-x}, \quad x\geq y
\end{cases}
$$
and the operator $H$ in the same Hilbert space with the integral 
kernel $H(x,y)=H_0(x,y)+e^{-x}e^{-y}$. 
Thus, $V=H-H_0$ is a rank one operator. In fact, $H_0$ and $H$ are
resolvents (with the spectral parameter $-1$) of the operator 
$-\frac{d^2}{dx^2}$ in $L^2(0,\infty)$ with the Dirichlet and Neumann
boundary conditions at zero. 

Krein shows that in this example $D(\lambda)$ is not a Hilbert-Schmidt
operator for $\lambda\in(0,1)$. A more detailed analysis \cite{KM} 
shows that the spectrum of $D(\lambda)$ is simple, purely a.c. and coincides with 
$[-1,1]$. 

What can be said about the scattering matrix in this case? 
First note that the spectra of both $H_0$ and $H$ are simple, purely a.c. and coincide with 
$[0,1]$. Thus, the fibre spaces $\h(\lambda)$ in \eqref{a12}
are one-dimensional and so the scattering matrix is simply a unimodular complex number. 
Krein calculates the spectral shift function $\xi(\lambda)$ for 
this pair of operators and shows that $\xi(\lambda)=1/2$ on $[0,1]$. 
Together with the Birman-Krein formula $\det S(\lambda)=e^{-2\pi i \xi(\lambda)}$ 
this shows that $S(\lambda)=-1$ for all $\lambda\in [0,1]$. Thus, we have 
a complete agreement with Theorems~\ref{th1} and \ref{th2}. 
It is not difficult to check that Hypotheses~\ref{hyp1} and \ref{hyp2} hold true with 
$\delta=(0,1)$.

\textbf{4. Example: Schr\"odinger operator.} 
Let $H_0=-\Delta$ in $\H=L^2(\R^d)$, $d=1,2,3$, and $H=H_0+V$, where $V$ is the operator of multiplication 
by a function $V:\R^d\to\R$, which is assumed to satisfy 
\begin{equation}
\abs{V(x)}\leq C(1+\abs{x})^{-\rho}, 
\quad \rho>1.
\label{a18}
\end{equation}
It is well known that under the assumption \eqref{a18}, the wave operators for the pair 
$H_0$ and $H$ exist and are complete, and the scattering matrix $S(\lambda)$ is well
defined and differs from the identity by a compact operator. 

\begin{theorem}
\label{th3}
(i) Assume $\rho>1$. Then for all $\lambda>0$, formula \eqref{a3} holds true. 

(ii) Assume $\rho>d$. Then for all $\lambda>0$, the conclusion of Theorem~\ref{th2} holds
true. 
\end{theorem}
To the best of the author's knowledge, this result is new even for  $d=1$.

\textbf{5. Acknowledgements.}
The author is grateful to D.~Yafaev for useful discussions. 
Part of the work was completed when the author stayed at the California Institute of Technology 
as a Leverhulme Fellow; the author is grateful to the Leverhulme Foundation for the financial 
support and to Caltech for hospitality.

\section{Proof of Theorem~\ref{th3}}

1. Fix $a<0$ such that $a<\inf \sigma(H)$. 
Consider the operators $h=(H-a)^{-1}$, $h_0=(H_0-a)^{-1}$. 
By the invariance principle for the scattering matrix (see \cite{Birman} or \cite{Yafaev}), 
we have
$$
S(\lambda; H,H_0)=S(\mu; h,h_0), 
\quad \mu=\frac{1}{\lambda-a},
\quad \lambda>0.
$$
Also, denoting by $E_{h_0}(\mu)$ and $E_h(\mu)$ the spectral projections of $h_0$ and $h$ 
associated with the interval $(-\infty, \mu)$, we have:
$$
E(\lambda)-E_0(\lambda)=E_{h_0}(\mu)-E_h(\mu), 
\quad \mu=\frac{1}{\lambda-a},
\quad \lambda>0.
$$
Thus, Theorem~\ref{th3} will follow from Theorems~\ref{th1} and \ref{th2} if we show that the pair of 
operators $h$, $h_0$ satisfies Hypothesis~\ref{hyp1} for $\rho>1$ and Hypothesis~\ref{hyp2} 
for $\rho>d$. 

In order to check this, we need to fix an appropriate factorization of $h-h_0$. We shall
use the factorization $h-h_0=g^* v_0 g$, where
$$
g=\abs{V}^{1/2}h_0, \quad 
v_0=-V_0-V_0 \abs{V}^{1/2} h \abs{V}^{1/2} V_0, \quad 
V_0=\sign(V). 
$$
This factorization is merely an iterated resolvent identity written in different notation. 

2. Assume $\rho>1$. It is well known that $\abs{V}^{1/2}h_0\in {\mathfrak S}_\infty$. 
We use the notation 
$$
T_0(z)=\abs{V}^{1/2}(H_0-z)^{-1}\abs{V}^{1/2}, 
\quad
T(z)=\abs{V}^{1/2}(H-z)^{-1}\abs{V}^{1/2},
\quad 
\Im z>0.
$$
By the spectral theorem, we have 
\begin{equation}
\frac{d}{d\mu} g E_{h_0}(\mu)g^*=(\lambda-a)^2\frac{d}{d\lambda}\abs{V}^{1/2} E_0(\lambda)\abs{V}^{1/2}
=(\lambda-a)^2 \frac{1}{\pi}\Im T_0(\lambda+i0).
\label{c4}
\end{equation}
The limit $T_0(\lambda+i0)$ exists and is continuous in $\lambda>0$ in the operator norm.
This fact is known as the limiting absorption principle; 
it stems from the Sobolev's embedding theorems. 

Next, we need to discuss the derivative $\frac{d}{d\mu}gE_h(\mu)g^*$. 
Before doing this, let us recall the following facts: 
\begin{gather}
T(z)=T_0(z)(I+V_0T_0(z))^{-1}, 
\quad \Im z>0,
\label{c1} 
\\
I+V_0T_0(\lambda+i0)\text{ has a bounded inverse for all $\lambda>0$. }
\label{c2}
\end{gather}
Formula \eqref{c1} follows from the resolvent identity. 
Relation \eqref{c2} goes back to Agmon \cite{Agmon} and uses the fact that (by Kato's theorem 
\cite{Kato}) $H$ has no positive eigenvalues. 
Also due to Agmon is the observation that one can put \eqref{c1} and \eqref{c2} together
and prove that $T(\lambda+i0)$ is H\"older continuous in $\lambda>0$ in the operator norm. 
It follows that 
\begin{equation}
\text{the derivative } 
\frac{d}{d\lambda}\abs{V}^{1/2} E(\lambda)\abs{V}^{1/2}
=\frac{1}{\pi}\Im T(\lambda+i0)
\text{ exists and is H\"older continuous.}
\label{b13}
\end{equation}

Let us return to the derivative $\frac{d}{d\mu}gE_h(\mu)g^*$. 
Using the resolvent identity, we get 
\begin{multline}
gE_h(\mu)g^*=
\abs{V}^{1/2}h_0 E_h(\mu) h_0\abs{V}^{1/2}
\\
=
\abs{V}^{1/2}h E_h(\mu) h\abs{V}^{1/2}
+
\abs{V}^{1/2}h_0Vh E_h(\mu) hVh_0\abs{V}^{1/2}
\\
+
\abs{V}^{1/2}h_0Vh E_h(\mu) h\abs{V}^{1/2}
+
\abs{V}^{1/2}h E_h(\mu) hVh_0\abs{V}^{1/2}.
\label{c3}
\end{multline}
Inspecting each term in the r.h.s. and using \eqref{b13}, 
we see that the derivative of the above 
expression exists and is H\"older continuous in $\mu$ in the operator norm.

3. Assume $\rho>d$. 
It is well known that $\abs{V}^{1/2}h_0\in{\mathfrak S}_2$; this follows from an inspection 
of the integral kernel of this operator. 

Next, we claim that the derivative $\frac{d}{d\lambda}\abs{V}^{1/2} E_0(\lambda)\abs{V}^{1/2}$
exists and is H\"older continuous in the trace norm. 
This fact is probably well known to specialists; in any case, it follows 
from a simple computation involving factorization of the pre-limiting expressions into products
of two Hilbert-Schmidt operators and estimating the Hilbert-Schmidt norm of each of these factors. 
The details of this computation can be found, e.g. in \cite{Push}. 
By \eqref{c4}, the derivative $\frac{d}{d\mu} g E_{h_0}(\mu)g^*$ also exists and is H\"older continuous
in the trace norm. 

Finally, consider the derivative $\frac{d}{d\mu}gE_h(\mu)g^*$. First, by using \eqref{c3} we reduce 
the question to the existence and H\"older continuity of 
$\frac{d}{d\lambda}\abs{V}^{1/2} E(\lambda)\abs{V}^{1/2}$.
The latter fact again follows from \eqref{c1}, \eqref{c2} and the H\"older continuity 
of $\frac{d}{d\lambda}\abs{V}^{1/2} E_0(\lambda)\abs{V}^{1/2}$.

\section{Proof of Theorems~\ref{th1} and \ref{th2}}
We use the notation
$$
R_0(z)=(H_0-zI)^{-1}, \quad R(z)=(H-zI)^{-1},\quad
T_0(z)=GR_0(z)G^*, \quad T(z)=GR(z)G^*.
$$
For $\lambda\in \delta$, let us introduce an auxiliary operator in $\K$: 
\begin{equation}
A(\lambda)=\pi^2(F'_0(\lambda))^{1/2}V_0 F'(\lambda) V_0 (F'_0(\lambda))^{1/2}.
\label{a7}
\end{equation}
Clearly, $A(\lambda)$ is compact, self-adjoint, and $A(\lambda)\geq0$.
This operator plays an important role in our construction. 
As we shall see later, the spectrum of $A(\lambda)$ is related to the spectrum of the scattering 
matrix $S(\lambda)$. In order to describe this relation, let us introduce the following  notation. 
For  bounded normal operators $X$ and $Y$ in Hilbert spaces 
$\H_X$ and $\H_Y$, we shall write 
$$
X\approx Y \text{ if } X\mid_{\H_X\ominus \Ker X}
\text{ is unitarily equivalent to }Y\mid_{\H_Y\ominus \Ker Y}.
$$
It is well known that $X^*X\approx XX^*$ for any bounded operator $X$;
we shall repeatedly use this fact.
\begin{lemma}\label{l3}
Suppose that the Hypothesis~\ref{hyp1} holds true. 
Then for all $\lambda\in\delta$,  
\begin{equation}
A(\lambda)\approx \frac14 (S(\lambda)-I_\lambda)^*(S(\lambda)-I_\lambda)=\frac12(I_\lambda-\Re S(\lambda)).
\label{a9}
\end{equation}
\end{lemma}
In other words, the Lemma says that if $e^{i\theta_n}$ are the eigenvalues of $S(\lambda)$, then 
 $(\sin(\theta_n/2))^2$ are the eigenvalues of $A(\lambda)$.
\begin{proof}
1. 
First we recall the stationary representation for the scattering matrix. 
For $f\in\H^{(ac)}_0(\d)$, let $\{f(\lambda)\}_{\lambda\in\d}$, $f(\lambda)\in\h(\lambda)$, 
be the representation of $f$ in the direct integral \eqref{a12}.
Then for all $\lambda\in\d$, the operator $\F(\lambda):\K\to \h(\lambda)$, 
$f\mapsto (G^* f)(\lambda)$ is well defined, bounded, and 
\begin{equation}
\F(\lambda)^*\F(\lambda)=F_0'(\lambda).
\label{a14}
\end{equation}
For a.e. $\lambda\in\delta$, the scattering matrix can be represented as
\begin{equation}
S(\lambda)=I_\lambda-2\pi i\F(\lambda)(V_0-V_0T(\lambda+i0)V_0)\F(\lambda)^*.
\label{a15}
\end{equation}

2. 
Consider an auxiliary unitary operator $\wt S(\lambda)$ 
in $\H$, defined by 
\begin{equation}
\wt S(\lambda)=I-2\pi i(F_0'(\lambda))^{1/2}(V_0-V_0T(\lambda+i0)V_0)(F_0'(\lambda))^{1/2}.
\label{a16}
\end{equation}
By virtue of \eqref{a14}, we have 
$$
S(\lambda)-I_\lambda\approx \wt S(\lambda)-I,
\quad \lambda\in\delta
$$
(see \cite[Lemma 7.7.1]{Yafaev}).
It follows that
\begin{equation}
(S(\lambda)-I_\lambda)^*(S(\lambda)-I_\lambda)\approx 
(\wt S(\lambda)-I)^*(\wt S(\lambda)-I).
\label{a17}
\end{equation}

3. 
For any $\e>0$,  employing the resolvent identity,  we obtain

\begin{multline*}
(V_0-V_0T(\lambda-i\e)V_0)(\Im T_0(\lambda+i\e))(V_0-V_0T(\lambda+i\e)V_0)
\\
=
V_0 G(I-R(\lambda-i\e)V)\e R_0(\lambda-i\e)R_0(\lambda+i\e)(I-VR(\lambda+i\e))G^*V_0
\\
=V_0G\e R(\lambda-i\e)R(\lambda+i\e) G^* V_0=V_0(\Im T(\lambda+i\e))V_0.
\end{multline*}
Taking $\e\to+0$ in the above identity and multiplying on 
both sides by $(F'_0(\lambda))^{1/2}$, we obtain
$$
\frac1{4\pi} (\wt S(\lambda)-I)^*(\wt S(\lambda)-I)
=
(F'_0(\lambda))^{1/2}V_0 (\Im T(\lambda+i0))^{1/2} V_0 (F'_0(\lambda))^{1/2}=\frac1\pi A(\lambda).
$$
Together with \eqref{a17}, this proves the required statement. 
\end{proof}

Let us fix $\lambda_0\in\delta$ and prove the conclusions of Theorems~\ref{th1} and \ref{th2} for this 
value $\lambda=\lambda_0$. In order to simplify our notation, let us assume (without the loss of generality) 
that $\lambda_0=0$.

We use the notation $\R_+=(0,\infty)$, $\R_-=(-\infty,0)$.

The proofs of Theorems~\ref{th1} and \ref{th2} will be deduced from
the following Lemma, which might be of some interest in its own right. 

\begin{lemma}\label{l4}
(i) Assume Hypothesis~\ref{hyp1} and $0\in\delta$. 
Then the essential spectra of the operators $E_0(\R_\pm)E(\R_\mp)E_0(\R_\pm)$ coincide
with $[0,\norm{A(0)}]$.

(ii) Assume Hypothesis~\ref{hyp2} and $0\in\delta$. 
Let $s_n$ be the non-zero eigenvalues of $A(0)$.
Then the a.c. parts of the operators $E_0(\R_\pm)E(\R_\mp)E_0(\R_\pm)$
are unitarily equivalent to a direct sum of operators of multiplication by $x$
in $L^2([0,s_n],dx)$.
\end{lemma}

\begin{proof}[Proof of Theorems \ref{th1} and \ref{th2}]

1. First let us reduce our considerations to the case
\begin{equation}
E_0(\{0\})=E(\{0\})=0.
\label{b1}
\end{equation}
Hypothesis~\ref{hyp1} for $\lambda=0$
implies that $GE_0(\{0\})G^*=0$; therefore, $GE_0(\{0\})=0$ and so 
$VE_0(\{0\})=0$. It follows that the subspace $E_0(\{0\})$ reduces both 
$H$ and $H_0$, and so $E_0(\{0\})=E(\{0\})$. Denote
$\wt \H=\H\ominus E_0(\{0\})$, $\wt H=H\mid_{\wt \H}$,  $\wt H_0=H_0\mid_{\wt \H}$,
and let $\wt D(0)$ be the difference \eqref{a1} constructed for the operators $\wt H_0$, $\wt H$.
Then we have $\wt D(0)\approx D(0)$ and $0$ is not an eigenvalue 
of $\wt H$ or of $\wt H_0$. Thus, without the loss of generality we can assume that 
from the start \eqref{b1} holds true.

2.
Let us denote $D=D(0)$ and 
\begin{align*}
\H_+=&\Ker(D-I)=\Ran E(\R_-)\cap \Ker E_0(\R_-),
\\
\H_-=&\Ker(D+I)=\Ran E_0(\R_-)\cap \Ker E(\R_-),
\\
\H_0=&\H\ominus(\H_-\oplus\H_+).
\end{align*}
It is well known (see \cite{Halmos} or \cite{ASS}) that 
$D\mid_{\H_0}\approx (-D)\mid_{\H_0}$.
Therefore, the spectral analysis of $D$ reduces to the spectral analysis
of $D^2$ and to the spectral analysis of the dimensions of $\H_+$ and 
$\H_-$. 
Next,  $D^2$ can be represented as  
\begin{equation}
D^2=E_0(\R_-)E(\R_+)E_0(\R_-)+E_0(\R_+)E(\R_-)E_0(\R_+),
\label{b11}
\end{equation}
and the r.h.s. provides a block-diagonal decomposition of $D^2$
with respect to the decomposition $\H=\Ran E_0(\R_-)\oplus \Ran E_0(\R_+)$.
Thus, the spectral analysis of $D^2$ reduces to the spectral analysis of the 
two terms on the r.h.s. of \eqref{b11}. 

3. 
Taking into account the decomposition \eqref{b11}, we see that 
Theorem~\ref{th2} follows directly from Lemma~\ref{l4}(ii). 

Similarly, Lemma~\ref{l4}(i) characterises $\sigma_{ess}(D)$ away 
from $-1$ and $1$. 
In order to complete the proof of Theorem~\ref{th1}, it remains to take care of the eigenvalues 
$\pm1$ of $D$. If $\norm{A(0)}<1$, then Lemma~\ref{l4}(i)  ensures that the 
kernels 
$$
\Ker (E_0(\R_\pm)E(\R_\mp)E_0(\R_\pm)-I)=\H_\pm
$$
are finite dimensional, and so $\pm1$ do not contribute to the essential spectrum
of $D$. On the other hand, if $\norm{A(0)}=1$, then 
by Lemma~\ref{l4}, $\s_{ess}(D)=[-1,1]$ regardless of  the 
dimensions of $\H_\pm$ and so we have nothing to prove.
\end{proof}

The key element in our proof of Lemma~\ref{l4} is a representation of the product
$E(\R_-)E_0(\R_+)$ in terms of some auxiliary operators $Z$, $Z_0$ which we proceed 
to define. These operators act from $L^2(\R_+,\K)$ into $\H$; here 
$L^2(\R_+,\K)$ is the space of  measurable functions  
$f:\R_+\to\K$ such that 
$$
\int_0^\infty \norm{f(t)}^2_{\K}dt<\infty.
$$
$L^1(\R_+,\K)$ is defined similarly. 
On the dense subset $L^1(\R_+,\K)\cap L^2(\R,\K)$, let us define the operators 
$Z$, $Z_0$ by 
\begin{align*}
Z_0 f=&\int_0^\infty e^{-tH_0} E_0(\R_+)G^*f(t)dt,
\\
Z f=&\int_0^\infty e^{tH} E(\R_-)G^*f(t)dt.
\end{align*}
We will see (in Lemma~\ref{l8}) that $Z_0$ and $Z$ are bounded and
\begin{equation}
E(\R_-)E_0(\R_+)=-ZV_0Z_0^*.
\label{b3}
\end{equation}
From \eqref{b3} we get the representation formula
$E_0(\R_+)E(\R_-)E_0(\R_+)=
Z_0 V_0 Z^* Z V_0 Z_0^*$, which will be important in our proof of Lemma~\ref{l4}. 
But first we need to develop some analysis related to the operators $Z$ and $Z_0$; 
this is done in the beginning of the next section.

\section{Hankel operators; Proof of Lemma~\ref{l4}}

\textbf{1. Hankel operators.}
We need to prepare some estimates for vector valued Hankel operators.  
These are straightforward generalisations of the well known technique 
of spectral theory of Hankel operators (see \cite{Power,Howland1,Howland2,Howland3})
to a vector valued case.

Suppose that for each $t\geq0$, a bounded self-adjoint operator $K(t)$ in $\K$ is given. 
Suppose that $K(t)$ is continuous in $t\geq0$ in the operator norm.
Define a Hankel type operator $K$ in $L^2(\R_+,\K)$ by
\begin{equation}
(Kf,g)_{L^2(\R_+,\K)}=\int_0^\infty \int_0^\infty 
(K(t+s)f(s),g(t))_{\K}dt\,ds,
\label{b6}
\end{equation}
when $f,g\in L^2(\R_+,\K)\cap L^1(\R_+,\K)$.
\begin{lemma}\label{l6}
(i) Suppose $\norm{K(t)}\leq C_1/t$ for all $t>0$.
Then the operator $K$ is bounded and $\norm{K}\leq \pi C$.
(ii) Suppose $K(t)$ is compact for all $t$ and $\norm{K(t)}t\to0$ as $t\to+0$ and 
as $t\to+\infty$. Then $K$ is compact.
(iii) Suppose 
$$
K(t)=\int_0^\infty M(\lambda) e^{-\lambda t} d\lambda,
$$
where $M(\lambda)$ is a measurable function of $\lambda\in(0,\infty)$ with values
in the set of trace class operators in $\K$. 
Suppose that 
$$
C_2:=\int_0^\infty \norm{M(\lambda)}_{{\mathfrak S}_1}\lambda^{-1}d\lambda<\infty; 
$$
then $K$ is a trace class operator. 
\end{lemma}
\begin{proof}
(i), (ii) is a straightforward generalisation of Proposition~1.1 from \cite{Howland3}.
Indeed, since the Carlemann operator on $L^2(\R_+)$ with the kernel $(t+s)^{-1}$ 
is bounded with the norm $\pi$, we have 
$$
\abs{(Kf,g)_{L^2}}\leq C_1 \int_0^\infty \int_0^\infty 
\frac{\norm{f(s)}_{\K}\norm{g(t)}_{\K}}{t+s}dt\,ds
\leq
\pi C_1 \norm{f}_{L^2(\R_+,\K)} \norm{g}_{L^2(\R_+,\K)},
$$
which proves (i). To prove (ii), we need to approximate $K$ by compact
operators. Let $K_n(t)=K(t)\chi_{(1/n,n)}(t)$ and let $K_n$ be the corresponding
operator in $L^2(\R_+,\K)$. It is not difficult to see that each $K_n$ is compact.
By (i), $\norm{K-K_n}\to0$ as $n\to\infty$.

(iii) For each $\lambda$, let us represent $M(\lambda)$ as a difference of its positive
and negative parts: $M(\lambda)=M_+(\lambda)-M_-(\lambda)$, 
$M_\pm(\lambda)\geq0$, 
$\norm{M(\lambda)}_{{\mathfrak S}_1}=\Tr M_+(\lambda)+\Tr M_-(\lambda)$.
Then $K$ splits accordingly as $K=K_+-K_-$.
Let us factorize each of $K_+$, $K_-$ into a product of Hilbert-Schmidt operators as follows. 
Let 
\begin{gather*}
N_\pm: L^2(\R_+,\K)\to L^2(\R_+,\K), 
\\
(N_\pm f)(\lambda)=M_\pm(\lambda)^{1/2} \int_0^\infty e^{-\lambda t} f(t) dt.
\end{gather*}
Then $K_\pm=N_\pm^* N_\pm$ and 
\begin{multline*}
\norm{N_\pm}^2_{{\mathfrak S}_2}=\int_0^\infty d\lambda \int_0^\infty dt \norm{e^{-\lambda t}M_\pm(\lambda)^{1/2}}_{{\mathfrak S}_2}^2
\\
=
\int_0^\infty d\lambda \norm{M_\pm(\lambda)}_{{\mathfrak S}_1} \int_0^\infty  e^{-2\lambda t}dt
=
\int_0^\infty \Tr (M_\pm(\lambda))(2\lambda)^{-1}d\lambda<\infty,
\end{multline*}
which yields the required result. 
\end{proof}

Consider the self-adjoint operators $\G_0$, $\G$ in $L^2(\R_+)$ 
which are given by the integral kernels 
$$
\G_0(t,s)=\frac{e^{-t-s}}{t+s},
\quad
\G(t,s)=\frac{1-e^{-t-s}}{t+s}.
$$
It is well known that $\G_0$ is bounded and has purely a.c. spectrum $[0,\pi]$
of multiplicity one; explicit diagonalisation of $\G_0$ is available (see \cite{Rosenblum}).
The following proposition is probably well known to specialists, but we were 
unable to find it in the literature. 
\begin{lemma}\label{l7}
The operator $\G$ is unitarily equivalent to $\G_0$.
Thus, $\Gamma$ has a purely a.c. spectrum of multiplicity one 
which coincides with $[0,\pi]$. 
\end{lemma}
\begin{proof}
The proof is a combination of identities from \cite{Howland2}.
Let $N: L^2(\R_+)\to L^2(\R_+)$ be the operator
$(Nf)(t)=\int_0^\infty e^{-ts}f(s)ds$.
We have $\G=N\chi_{(0,1)}N$, $\G_0=N\chi_{(1,\infty)}N$.
Next, let $U: L^2(\R_+)\to L^2(\R_+)$ be the unitary operator 
$(Uf)(x)=\frac1x f(1/x)$.
Then $U^2=I$, $U\chi_{(0,1)}=\chi_{(1,\infty)}U$ and $UN^2 U=N^2$.
Using the well known fact that $X^*X\approx X X^*$, we get 
\begin{multline*}
\G=N\chi_{(0,1)}N=(N\chi_{(0,1)}U)(U\chi_{(0,1)}N)=
(NU\chi_{(1,\infty)})(\chi_{(1,\infty)}UN)
\\
\approx (\chi_{(1,\infty)}UN)(NU\chi_{(1,\infty)})
=
(\chi_{(1,\infty)}N)(N\chi_{(1,\infty)})
 \approx N\chi_{(1,\infty)}N=\G_0.
 \end{multline*}
 Thus, $\Gamma\approx\Gamma_0$.
 It remains to note that $\Ker \G=\Ker \G_0=\{0\}$.
 \end{proof}

 \textbf{2. Proof of Lemma~\ref{l4}.}
 Important ``model'' operators in our considerations are the integral 
 Hankel operators in $L^2(\R_+,\K)$ of the type \eqref{b6} with the kernels
 given by 
 $$
 K(t)=F'(0)\frac{1-e^{-t}}{t}, 
 \quad 
 K_0(t)=F_0'(0)\frac{1-e^{-t}}{t}.
 $$
 Identifying $L^2(\R_+,\K)$ with $L^2(\R_+)\otimes \K$, we will denote these 
 operators by $\G\otimes F'(0)$ and $\G\otimes F_0'(0)$.
 
 \begin{lemma}\label{l8}
 (i) Assume that 
 \begin{equation}
\norm{F_0(\lambda)-F_0(0)}=O(\lambda),
\quad
\norm{F(\lambda)-F(0)}=O(\lambda),
\quad 
\text{as $\lambda\to0$.}
\label{b2}
\end{equation}
Then the operators $Z$ and $Z_0$ are bounded.
 
 (ii) Assume Hypothesis~\ref{hyp1} with $0\in\delta$. 
  Then the differences
 \begin{equation}
 Z_0^* Z_0-(\G\otimes F_0'(0))
 \quad \text{and} \quad
 Z^* Z-(\G\otimes F'(0))
 \label{b7} 
 \end{equation}
 are compact. 
 
 (iii) Assume Hypothesis~\ref{hyp2} with $0\in\delta$. 
 Then the differences \eqref{b7} are trace class operators.
 \end{lemma}
\begin{proof}
We will prove the statements for $Z_0$; the proofs for $Z$ are analogous.

(i) Let $f\in L^2(\R_+,\K)\cap L^1(\R_+,\K)$; we have 
$$
\norm{Z_0 f}^2
=
\int_0^\infty \int_0^\infty
(Ge^{-(t+s)H_0}E_0(\R_+)G^*f(t),g(s))_{\K}dt\,ds,
$$
and so the above expression is a quadratic form of the operator of the type 
\eqref{b6} with the kernel 
$K(t)=Ge^{-tH_0}E_0(\R_+)G^*$.
By Lemma~\ref{l6}, it suffices to prove the bound 
$\norm{K(t)}\leq C/t$, $t>0$.
Using our assumption \eqref{b2}, we have 
\begin{multline*}
\norm{K(t)}=\nnorm{t\int_0^\infty e^{-t\lambda}G E_0((0,\lambda))G^* d\lambda}
\leq
t\int_0^\infty e^{-t\lambda}\norm{F_0(\lambda)-F_0(0)}d\lambda
\\
\leq Ct \int_0^\infty e^{-t\lambda}\lambda d\lambda=C/t.
\end{multline*}
(ii) By the same reasoning, $Z_0^*Z_0-(\G\otimes F_0'(0))$ is an operator of the type
\eqref{b6} with 
$$
K(t)=Ge^{-tH_0} E_0(\R_+) G^*-F_0'(0)\int_0^1 e^{-t\lambda}d\lambda.
$$
By Lemma~\ref{l6}, it suffices to prove that $\norm{K(t)}t\to0$ as $t\to0$ and $t\to\infty$.
For $t\to0$ this is clearly true. Next, we have 
\begin{equation}
K(t)=t\int_0^\infty e^{-t\lambda} (F_0(\lambda)-F_0(0))d\lambda
-F_0'(0)t\int_0^\infty \min\{\lambda,1\}e^{-t\lambda}d\lambda
\label{b8}
\end{equation}
and from $\norm{F_0(\lambda)-F_0(0)-\lambda F_0'(0)}=o(\lambda)$,
$\lambda\to0$, we conclude that $\norm{K(t)}=o(1/t)$ as $t\to\infty$.

(iii) 
Choose $\gamma>0$ such that $[0,\gamma]\subset\delta$. 
As above, $Z_0^*Z_0-(\G\otimes F_0'(0))$ has the kernel \eqref{b8}.
Let us write this kernel as $K(t)=K_1(t)+K_2(t)$, with 
\begin{multline*}
K_1(t)=G e^{-tH_0} E_0([0,\gamma]) G^*-F'_0(0)\int_0^1 e^{-t\lambda} d\lambda
\\
=
\int_0^\infty [F'_0(\lambda)\chi_{(0,\gamma)}(\lambda)-F'_0(0)\chi_{(0,1)}(\lambda)]e^{-t\lambda}d\lambda,
\end{multline*}
and 
$K_2(t)=Ge^{-t H_0} E_0((\gamma,\infty))G^*$, so $K_2=(E_0((\gamma,\infty))Z_0)^*(E_0((\gamma,\infty))Z_0)$.
By the H\"older continuity assumption, 
$$
\int_0^\infty
\norm{F'_0(\lambda)\chi_{(0,\gamma)}(\lambda)-F'_0(0)\chi_{(0,1)}(\lambda)}_{{\mathfrak S}_1}\lambda^{-1}d\lambda<\infty,
$$
and so by Lemma~\ref{l6}(iii), the Hankel operator with the kernel $K_1$ is trace class. 
Finally, it is easy to see that $E_0((\gamma,\infty))Z_0\in{\mathfrak S}_2$, since 
$$
\norm{E_0((\gamma,\infty))Z_0}_{{\mathfrak S}_2}^2=\int_0^\infty dt \norm{e^{-tH_0} E_0((\gamma,\infty))G^*}_{{\mathfrak S}_2}^2
\leq \int_0^\infty e^{-2\gamma t}\norm{G^*}_{{\mathfrak S}_2}^2 dt<\infty,
$$
and so the Hankel operator with the kernel $K_2$ also belongs to the trace class. 
This argument borrows its main idea from  \cite{Howland1}.
\end{proof}

\begin{lemma}\label{l5}
Assume \eqref{b2} and $E_0(\{0\})=E_0(\{0\})=\{0\}$. 
Then the identity \eqref{b3} holds true. 
\end{lemma}
\begin{proof}
Let $\gamma>0$ and let $\psi,\psi_0\in\H$ be vectors such that 
$E((-\g,0))\psi=E_0((0,\g))\psi_0=0$. Since the set of such vectors is dense in $\H$, 
it suffices to prove that 
\begin{equation}
(E_0(\R_+)\psi_0,E(\R_-)\psi)=-(V_0Z_0^*\psi_0,Z^*\psi)_{L^2(\R_+,\K)}
\label{b4}
\end{equation}
for all such vectors $\psi, \psi_0$.
For $\psi$ and $\psi_0$ of this class, $Z_0^*\psi_0$ and $Z^*\psi$ are given by 
\begin{align*}
(Z_0^*\psi_0)(t)&= Ge^{-tH_0}E_0(\R_+)\psi_0,
\\
(Z^*\psi)(t)&=Ge^{tH}E(\R_-)\psi,
\end{align*}
and so we have 
\begin{multline}
(V_0Z_0^*\psi_0,Z^*\psi)_{L^2(\R_+,\K)}=
\int_0^\infty (V_0 Ge^{-tH_0}E_0(\R_+)\psi_0,Ge^{tH}E(\R_-)\psi)_{\K}dt
\\
=
\int_0^\infty (Ve^{-tH_0}E_0(\R_+)\psi_0,e^{tH}E(\R_-)\psi)dt.
\label{b5}
\end{multline}
Consider the function $L(t)=(e^{-tH_0}E_0(\R_+)\psi_0,e^{tH}E(\R_-)\psi)$.
This function is continuous in $t\geq0$ and we have 
$L(0)=(E_0(\R_+)\psi_0,E(\R_-)\psi)$, $L(+\infty)=0$, 
and $L'(t)=(Ve^{-tH_0}E_0(\R_+)\psi_0,e^{tH}E(\R_-)\psi)$.
Combining this with \eqref{b5}, we get \eqref{b4}.
\end{proof}

\begin{proof}[Proof of Lemma~\ref{l4}]
We will prove the statement for $E_0(\R_+)E(\R_-)E_0(\R_+)$; the proof
for $E_0(\R_-)E(\R_+)E_0(\R_-)$ is analogous.
By Lemma~\ref{l5}, we have 
$$
E_0(\R_+)E(\R_-)E_0(\R_+)
=
Z_0 V_0 Z^* Z V_0 Z_0^*.
$$
Next, by Lemma~\ref{l8}(ii), for some compact operators $X_0$ and $X$ we have
\begin{equation}
Z_0 V_0 Z^* Z V_0 Z_0^*=Z_0 V_0 (\G\otimes F'(0))  V_0 Z_0^* +X,
\label{b9}
\end{equation}
\begin{multline}
Z_0 V_0 (\G\otimes F'(0))  V_0 Z_0^*
=
\left\{Z_0 V_0 (\G^{1/2}\otimes (F'(0))^{1/2})\right\}
\left\{ (\G^{1/2}\otimes (F'(0))^{1/2})  V_0 Z_0^*\right\}
\\
\approx 
\left\{ (\G^{1/2}\otimes (F'(0))^{1/2})  V_0 Z_0^*\right\}
\left\{Z_0 V_0 (\G^{1/2}\otimes (F'(0))^{1/2})\right\}
\\
=
 (\G^{1/2}\otimes (F'(0))^{1/2})  V_0 (\G\otimes F'_0(0))    V_0 (\G^{1/2}\otimes (F'(0))^{1/2})
+X_0
\\
=
\pi^{-2} \G^2\otimes A(0)+X_0.
\label{b10}
\end{multline}
Thus, by Weyl's theorem on the stability of the essential spectrum under the 
compact perturbations, 
$$
\s_{ess}(E_0(\R_+)E(\R_-)E_0(\R_+))\setminus\{0\}
=
\s_{ess}(\pi^{-2} \G^2\otimes A(0))\setminus\{0\}.
$$
By Lemma~\ref{l7},  the essential spectrum of $\pi^{-2} \G^2\otimes A(0)$ coincides with 
$[0,\norm{A(0)}]$. This  proves part (i) of the Lemma.

Next, assuming Hypothesis~\ref{hyp2} and using part (iii) instead of part (ii) 
of Lemma~\ref{l8},  we arrive at \eqref{b9}, \eqref{b10}
with $X_0$ and $X$ of the trace class. 
Thus, by the Kato-Rosenblum theorem on the stability of the a.c. spectrum 
under trace class perturbations, 
the a.c. part of $Z_0 V_0 Z^* Z V_0 Z_0^*$ is unitarily equivalent to the a.c. part
of $\pi^{-2}\G^2\otimes A(0)$. 
By Lemma~\ref{l7},
the latter operator is unitarily equivalent to a direct sum 
of operators of multiplication by $x$ in $L^2([0,s_n],dx)$, where $s_n$ are the 
eigenvalues of $A(0)$.
\end{proof}

\end{document}